\documentclass[10pt]{article}
\usepackage[english]{babel} %

\usepackage{amssymb}
\usepackage{amsmath}
\usepackage{amsthm}
\usepackage{enumitem}
\usepackage{url}

\theoremstyle{plain}
\newtheorem{theorem}{Theorem}
\newtheorem*{theorem*}{Theorem} 
\newtheorem{example}[theorem]{Example}

\newtheorem{question}[theorem]{Question}

\newtheorem{lemma}[theorem]{Lemma}
\newtheorem{corollary}[theorem]{Corollary}

\theoremstyle{definition}
\newtheorem{remark}[theorem]{Remark}

\def\gr{\operatorname{gr}}

\title{Complete factorizations of finite groups}
\author{Mikhail Kabenyuk}
\date{}

\begin{document}\maketitle

\begin{abstract}
Let $G$ be a group.
The subsets $A_1,\ldots,A_k$ of $G$ form a complete factorization of group $G$ 
if they are pairwise disjoint and each element $g\in G$ is uniquely represented as $g=a_1\ldots a_k$ 
with $a_i\in A_i$.
We prove the following theorem.

\begin{theorem*}
    Let $G$ be a finite nilpotent group.
    If $|G|=m_1\ldots m_k$ where $m_1,\ldots,m_k$ are integers greater $1$ and $k>2$,
    then there exist subsets $A_1,\ldots,A_k$ of $G$ which form a complete factorization of group $G$
    and $|A_i|=m_i$ for all $i=1,2,\ldots,k$.
\end{theorem*}
In addition, we give several examples of building complete factorization for some groups and formulate one open question.
\end{abstract}

\section{Introduction}\label{section-Introduction}
\def\Hajos{Haj$\acute{\rm{o}}$s}
\def\Szabo{Szab$\acute{\rm{o}}$}

Let $(A_1,\ldots,A_k)$ be a $k$-tuple of subsets of a group $G$.
If each element $g\in G$ is uniquely represented as $g=a_1\ldots a_k$ with $a_i\in A_i$, 
then we say  $(A_1,\ldots,A_k)$ forms a \textit{$k$-factorization} of the group $G$.
If $(A_1,\ldots,A_k)$ form a $k$-factorization of $G$ and $|A_i|=m_i$, then $|G|=m_1\ldots m_k$.

The usefulness of factorizations of finite abelian groups was demonstrated by \Hajos \cite{Hajos} 
by proving a classical Minkowski's conjecture \cite{Minkowski}
on lattice tilings, which says that in any lattice tiling of space by cubes, there are two cubes that meet face to face.
Details of classical and new results on factorizations of abelian groups can be found in books 
by \Szabo\ \cite{Szabo} and by \Szabo\ and Sands \cite{Szabo-Sands}.

For non-Abelian groups, the conjecture that every finite group 
admits a factorization in which every factor is of simple order or of order 4 has not yet been proved.
Such factorizations are called a minimal logarithmic signature (MLS).
It is known that the
existence of MLS for any finite group can be reduced to
the existence of MLS for finite simple groups
\cite{Rahimipour}.

Call a finite group $G$ $k$-factorizable 
if for any factorization $|G|=m_1\cdots m_k$ with $m_i>1$
there exist subsets $A_i$ of $G$ with $|A_i|=m_i$ such that
$G=A_1\cdots A_k$.
In \cite{kabenyuk} it was proved that for every integer $k\geq3$
there exists a finite group $G$ such that $G$ is not $k$-factorizable.
The same paper gives a complete list of groups of order at most $100$ 
that are not $k$-factorizable for at least one $k$. 

Recently, in a paper \cite{Chin-Wang-Wong} the notion of complete factorization of a finite group was introduced.
The subsets $A_1,\ldots,A_k$ of a group $G$ form a complete factorization of $G$
if they are pairwise disjoint and each element $g\in G$ is uniquely represented as $g=a_1\ldots a_k$ with $a_i\in A_i$.
In [1], the authors use a curious and original notion of 'systems of product sets' to prove the following statement \cite[Theorem 2.9]{Chin-Wang-Wong}.
If $G$ is a finite abelian group and $|G|=m_1\ldots m_k$ 
where $m_1,\ldots,m_k$ are integers greater $1$ and $k>2$, then 
there exist subsets $A_1,\ldots,A_k$ of $G$ which form a complete factorization of group $G$
and $|A_i|=m_i$ for all $i=1,2,\ldots,k$.
Note that there does not exist a complete factorization of a finite abelian group 
consisting of two factors \cite[Lemma 2.7]{Chin-Wang-Wong}.

Here we extend this statement to nilpotent groups and prove the following theorem.    
\begin{theorem*}[On the complete factorization of finite nilpotent groups]
    Let $G$ be a finite nilpotent group.
    If $|G|=m_1\ldots m_k$ where $m_1,\ldots,m_k$ are integers greater $1$ and $k\geq3$,
    then there exist subsets $A_1,\ldots,A_k$ of $G$ which form a complete factorization of group $G$
    and $|A_i|=m_i$ for all $i=1,2,\ldots,k$.
\end{theorem*}
Our proof differs from the proof of the authors of paper \cite{Chin-Wang-Wong} and thus we obtain 
a new proof of the theorem on the complete factorization of finite abelian groups.
In addition, we give several examples of building complete factorization for some groups and formulate one open question.

\phantom{a}

If $X$ is a subset of a group $G$, $|X|$ will denote the order of $X$,
and $\gr(X)$ will denote the subgroup of $G$ generated by $X$.
If $X=\{x\}$, then $\gr(x)$ stands for cyclic subgroup of $G$ generated by $x$.
We write $H<G$ or $H\leq G$ whenever $H$ is a subgroup of $G$.
If $H<G$, then $|G:H|$ is the index of $H$ in $G$.

\section{Proof of the theorem}
So, let $G$ be a finite nilpotent group of order $n$ and $n=m_1\ldots m_k$,
where $m_1,\ldots,m_k$ are integers greater $1$ and $k>2$.
Since the converse of Lagrange's theorem is true for finite nilpotent groups, 
there exists a subgroup series
$$
H_1<H_2<\ldots<H_k=G
$$
such that $|H_i|=m_1\ldots m_i$ for all $i=1,2,\ldots,k$.
Let $T_i$ be a right transversal to $H_{i-1}$ in $H_i$ for all $i=2,\ldots,k$ and $T_1=H_1$.
We have $H_{i-1}T_i=H_i$ and
\begin{equation}\label{transversal}
T_1T_2T_3\ldots T_k=H_1T_2T_3\ldots T_k=H_2T_3\ldots T_k=\ldots=H_k=G.
\end{equation}
However, the sets $\{T_i\}$ do not form a complete factorization of the group $G$, since for example $T_2\cap T_1\neq\varnothing$.
We have to fix the sets $T_i$.

We will choose $h_i\in N_{H_{i+1}}(H_i)\setminus H_i$, where $N_{H_{i+1}}(H_i)$ is
the normalizer of $H_i$ in $H_{i+1}$.
Due to the normalizer property of nilpotent groups $H_i$ is a proper subgroup of $N_{H_{i+1}}(H_i)$
and hence the set $N_{H_{i+1}}(H_i)\setminus H_i$ is not empty.
Let 
$$
A_i=T_ih_i,\ i=2,3,\ldots,k-1 \text{ and } A_1=T_1.
$$
The sets $A_i$ have the following elementary properties
in the future we will use them not always explicitly referring to them:
\begin{enumerate}[label=(\roman*)]
  \item\label{properties_cardinal} 
  $|A_i|=m_i$, $i=1,2,\ldots,k-1$;
  \item\label{properties_Ai_subset_H_i+1}
  $A_i=T_ih_i\subset H_{i+1}$, $i=2,\ldots,k-1$;
  \item\label{properties_H_i-1Ai=Hihi} 
  $H_{i-1}A_i=H_{i-1}(T_ih_i)=(H_{i-1}T_i)h_i=H_ih_i$, $i=2,\ldots,k-1$;
  \item\label{properties_empty} 
  if $1\leq i<j\leq k-1$, then $A_i\cap A_j=\varnothing$.
\end{enumerate}
Let us prove the last property that 
$A_1,A_2,\ldots,A_{k-1}$ is pairwise disjoint.
Indeed, by virtue of properties \ref{properties_Ai_subset_H_i+1} and \ref{properties_H_i-1Ai=Hihi} we obtain
\begin{equation}\label{A_i_and_A_j}
    A_i\subset H_{i+1}\leq H_j\hbox{ and }
    A_j\subset H_jh_j.
\end{equation}
Since by choice $h_j\notin H_j$, it follows from $(\ref{A_i_and_A_j})$ that $A_i\cap A_j=\varnothing$.

We are left to construct $A_k$.
Let $t,s\in T_k$ be representatives of the cosets $H_{k-1}$ and $H_{k-1}h_{k-1}$, respectively.
In other words $\{t\}=T_k\cap H_{k-1}$ and $\{s\}=T_k\cap H_{k-1}h_{k-1}$.

Now pick 
$$
t'\in H_2\setminus H_1
$$
and 
$$
s'\in H_{k-1}h_{k-1}\setminus A_{k-1}.
$$
Note that this choice of $t'$ and $s'$ is possible since $|H_2|>|H_1|$ and for $k\geq3$ we have 
$$
m_1\ldots m_{k-1}>m_{k-1}\Rightarrow
|H_{k-1}h_{k-1}|>|A_{k-1}|.
$$
Therefore, $H_2\setminus H_1\neq\varnothing$ and $H_{k-1}h_{k-1}\setminus A_{k-1}\neq\varnothing$.

Let 
$$
A_k=\left(T_k\setminus\{t,s\}\right)\cup\{t',s'\}.
$$
It is clear that $A_k$ as well as $T_k$ is a set of right coset representatives of $H_{k-1}$ in $H_k$
and hence 
$$
A_k\cap H_{k-1}=\{t'\}\text{ and }A_k\cap H_{k-1}h_{k-1}=\{s'\}.
$$
We must check that $A_k\cap A_i=\varnothing$ for $i=1,2,\ldots,k-1$.
Let us consider 3 cases.

\textsc{Case 1.}
Let $i=1$.
We have $A_k\cap A_1\subset A_k\cap H_{k-1}=\{t'\}$.
If $A_k\cap A_1\neq\varnothing$, then $t'\in A_1=H_1$
that is contrary to the choice of $t'$.

\textsc{Case 2.}
If $2\leq i\leq k-2$, then since $A_i\subset H_ih_i\subset H_{i+1}\leq H_{k-1}$, 
we get 
$$
A_k\cap A_i\subset A_k\cap H_ih_i\subset A_k\cap H_{i+1}\subset A_k\cap H_{k-1}=\{t'\}.
$$
If $A_k\cap A_i\neq\varnothing$, then $t'\in A_i\subset H_ih_i$.
However, it contradicts that $t'\in H_2\leq H_i$ and $H_i\cap H_ih_i=\varnothing$.
It follows that $A_k\cap A_i=\varnothing$.

\textsc{Case 3.} It remains to prove that $A_k\cap A_{k-1}=\varnothing$.
Since $A_{k-1}\subset H_{k-1}h_{k-1}$ and $A_k$ is a set of right coset representatives of $H_{k-1}$ in $H_k$, it follows that 
$$
A_k\cap A_{k-1}\subset A_k\cap H_{k-1}h_{k-1}=\{s'\}.
$$
But due to the choice $s'\notin A_{k-1}$ it means 
$A_k\cap A_{k-1}=\varnothing$.

By virtue of property \ref{properties_cardinal}, $|A_i|=m_i$ for $i=1,2,\ldots,k-1$, and 
$|A_k|=|T_k|=|H_k:H_{k-1}|=m_k$.
Furthermore, using (\ref{transversal}) and 
the equality $h_i^{-1}H_ih_i=H_i$ or $H_ih_i=h_iH_i$ for all $i=2,\ldots,k-1$, 
we successively obtain:
\begin{align*}
A_1A_2A_3\ldots A_{k-1}A_k
&=
T_1T_2h_2T_3h_3\ldots T_{k-1}h_{k-1}A_k\\
&=
H_2h_2T_3h_3\ldots T_{k-1}h_{k-1}A_k\\
&=
h_2H_2T_3h_3\ldots T_{k-1}h_{k-1}A_k\\
&\ldots\\
&=
(h_2h_3\ldots h_{k-1})H_{k-1}A_k\\
&=
(h_2h_3\ldots h_{k-1})H_k\\
&=
G.
\end{align*}

This completes the proof of the theorem.

\section{Examples}
We will present some examples of complete factorizations as an illustration of the proof of the theorem.
Let us start with a case that is very simple.
\begin{example}
    Let $Q_8$ be the quaternion group.
    Construct a complete factorization $A_1,A_2,A_3$ of the group $Q_8$ with $|A_i|=2$, $i=1,2,3$.
\end{example}
The quaternion group $Q_8$ consists of eight quaternions: $\pm1,\pm i,\pm j, \pm k$.
Let $H_1=\langle-1\rangle$, $H_2=\langle i\rangle$, and $H_3=Q_8$ and 
$T_1=H_1$, $T_2=\{1,i\}$, and $T_3=\{1,j\}$.
If $h_2=j$, then $A_2=T_2h_2=\{j,k\}$.
To construct $A_3$ we must choose $t'\in H_2\setminus H_1=\{i,-i\}$ and $s'\in H_2h_2\setminus A_2=\{-j,-k\}$.
Let $t'=i$ and $s'=-j$. 
We get $A_3=\{i,-j\}$ and if $A_1=T_1=H_1=\{1,-1\}$, then 
\begin{align*}
A_1A_2A_3 
&=\{1,-1\}\cdot\{j,k\}\cdot\{i,-j\}\\
&=\{\pm j,\pm k\}\cdot\{i,-j\}=\{\pm k\pm j\}\cup\{\pm1,\pm i\}=Q_8.
\end{align*}

\begin{example}
    Let $G=\mathbb{Z}_2^n$, $n>2$.
    Construct a complete factorization $A_1,\ldots,A_n$ of the group $G$ that $|A_i|=2$, $i=1,\ldots,n$.
\end{example}
The elements of the group $G$ are binary vectors of length $n$ having the form $(\alpha_1,\ldots,\alpha_n)$ where $\alpha_i=0$ or $1$.
Let us denote by $e_i$ a binary vector
each of whose components are all zero, except $i$-th that equals $1$.
Let $H_i=\gr(e_1,\ldots,e_i)$ and $T_i=\{0,e_i\}$, $i=1,\ldots,n$. 
Then let $h_i=e_{i+1}$ and 
$$
A_i=T_i+h_i=\{e_{i+1},e_i+e_{i+1}\},\ i=2,3,\ldots,n-1. 
$$
Note that $T_n\cap H_{n-1}=\{0\}$ and $T_n\cap H_{n-1}+h_{n-1}=T_n\cap H_{n-1}+e_n=\{e_n\}$.
Hence $t=0$ and $s=e_n$.
Now choose $t'\in H_2\setminus H_1$ and $s'\in H_{n-1}+e_n\setminus A_{n-1}$
(by our assumption $n-1\neq1$). 
We can take, for instance, $s'=e_1+e_n$ and $t'=e_2$.
Now we have no other choice $A_n=\{e_2,e_1+e_n\}$.
Here is a complete list of all $A_i$ which form a complete factorization of the group $\mathbb{Z}_2^n$:
\begin{align*}
&A_1=\{0,e_1\},\\
&A_2=\{e_3,e_2+e_3\},\\
&\ldots\\
&A_{n-1}=\{e_n,e_{n-1}+e_n\},\\
&A_n=\{e_2,e_1+e_n\}.
\end{align*}

\begin{example}
    Let $G=\mathbb{Z}_{10^n}$ and $n\geq3$.    
    Construct a complete factorization $A_1,\ldots,A_n$ of the group $G$ that $|A_i|=10$, $i=1,\ldots,n$.
\end{example}
We will denote the elements of the group $G$ by the symbols: $0,1,\ldots,10^n-1$. 
Pick the subgroups 
$$
H_1=\gr(10^{n-1})<\ 
H_2=\gr(10^{n-2})<\ 
\ldots<\ 
H_{n-1}=\gr(10)<\ 
H_n=G
$$ 
and the transversals $T_i$ to $H_{i-1}$ in $H_i$:
$$
T_i=10^{n-i}\cdot\{0,1,2,\ldots,9\},\ i=1,2,\ldots,n.
$$
Every integer $a$ in the interval $[0,10^n-1]$ uniquely appears as $a=t_1+\ldots+t_n$, where $t_i\in T_i$. This is just a decimal notation of the number $a$.

Now construct $A_i$, $i=1,\ldots,n$. Let 
$$
h_i=10^{n-i-1}\in H_{i+1}\setminus H_i,\ i=2,3,\ldots,n-1.
$$
Let's define $A_i$ by the following rules:
\begin{align*}
&A_1=T_1,\\
&A_i=T_i+10^{n-i-1}=10^{n-i-1}\cdot\{1,11,\ldots,91\},\ i=2,\ldots,n-1.
\end{align*}
To construct $A_n$, we choose $t$ and $s$ such that $\{t\}=T_n\cap H_{n-1}$ and $\{s\}=T_n\cap H_{n-1}+h_{n-1}$. 
We see that $t=0$ and $s=1$. 
Now let us choose $t'\in H_2\setminus H_1$ and $s'\in H_{n-1}+h_{n-1}\setminus A_{n-1}$ say $t'=10^{n-2}$, $s'=101$.
Therefore, we see that the following sets form a complete factorization of the group $\mathbb{Z}_{10^n}$:
\begin{align*}
&A_1=10^{n-1}\cdot\{0,1,2,\ldots,9\},\\
&A_2=10^{n-3}\cdot\{1,11,21,\ldots,91\},\\
&\ldots\\
&A_{n-2}=10\cdot\{1,11,21,\ldots,91\},\\
&A_{n-1}=\{1,11,21,\ldots,91\},\\
&A_n=\{10^{n-2},101,2,\ldots,9\}.
\end{align*}

\begin{example}
    Let $G=UT_3(\mathbb{F}_p)$ be the group of upper unitriangular matrices over a field $\mathbb{F}_p$,
    where $p$ is a prime.
    Construct a complete factorization $A_1,A_2,A_3$ of the group $G$ with $|A_i|=p$, $i=1,2,3$.
\end{example}
Let
$$
H_1=
\left(
  \begin{array}{ccc}
    1 & 0 & * \\
      & 1 & 0 \\
      &   & 1 \\
  \end{array}
\right),
H_2=
\left(
  \begin{array}{ccc}
    1 & 0 & * \\
      & 1 & * \\
      &   & 1 \\
  \end{array}
\right),
H_3=
\left(
  \begin{array}{ccc}
    1 & * & * \\
      & 1 & * \\
      &   & 1 \\
  \end{array}
\right),
$$
and
$$
T_1=
\left(
  \begin{array}{ccc}
    1 & 0 & * \\
      & 1 & 0 \\
      &   & 1 \\
  \end{array}
\right),
T_2=
\left(
  \begin{array}{ccc}
    1 & 0 & 0 \\
      & 1 & * \\
      &   & 1 \\
  \end{array}
\right),
T_3=
\left(
  \begin{array}{ccc}
    1 & * & 0 \\
      & 1 & 0 \\
      &   & 1 \\
  \end{array}
\right).
$$
Clearly, $H_3=G$ and $H_1,H_2$ are normal subgroups in $G$ and 
$T_1=H_1$, $H_1T_2=H_2$, and $H_2T_3=H_3$.
Next let 
$$
h_2=
\left(
  \begin{array}{ccc}
    1 & 1 & 0 \\
      & 1 & 0 \\
      &   & 1 \\
  \end{array}
\right).
$$ 
Now we will put $A_1=H_1$, $A_2=T_2h_2$. 
We choose $t'\in H_2\setminus H_1$ and $s'\in H_2h_2\setminus A_2$:
$$
t'=
\left(
  \begin{array}{ccc}
    1 & 0 & 0 \\
      & 1 & 1 \\
      &   & 1 \\
  \end{array}
\right),\
s'=
\left(
  \begin{array}{ccc}
    1 & 1 & 1 \\
      & 1 & 0 \\
      &   & 1 \\
  \end{array}
\right).
$$
So, we obtain
$$
A_1=
\left(
  \begin{array}{ccc}
    1 & 0 & * \\
      & 1 & 0 \\
      &   & 1 \\
  \end{array}
\right),
A_2=
\left(
  \begin{array}{ccc}
    1 & 1 & 0 \\
      & 1 & * \\
      &   & 1 \\
  \end{array}
\right),
$$
and
$$
A_3=
\left\{
\left(
  \begin{array}{ccc}
    1 & \alpha & 0 \\
      & 1 & 0 \\
      &   & 1 \\
  \end{array}
\right),\alpha\neq0,1,\ 
\left(
  \begin{array}{ccc}
    1 & 0 & 0 \\
      & 1 & 1 \\
      &   & 1 \\
  \end{array}
\right),\
\left(
  \begin{array}{ccc}
    1 & 1 & 1 \\
      & 1 & 0 \\
      &   & 1 \\
  \end{array}
\right)
\right\}.
$$

\begin{remark}
  The following formula illustrates the rule for using asterisks:
$$
H_2=
\left(
  \begin{array}{ccc}
    1 & * & * \\
      & 1 & 0 \\
      &   & 1 \\
  \end{array}
\right)=
\left\{\left(
  \begin{array}{ccc}
    1 & \alpha & \beta \\
      & 1 & 0 \\
      &   & 1 \\
  \end{array}
\right)\mid \alpha, \beta\in\mathbb{F}_p
\right\}.
$$  
Also, the empty spaces under the main diagonal here and in the example below must be filled with zeros.  
\end{remark}

\begin{example}
    Let $G=UT_n(\mathbb{F}_q)$ with $n>3$ be the group of upper unitriangular matrices over a field $\mathbb{F}_q$,
    where $q$ is a prime power.
    Construct a complete factorization $A_1,\ldots,A_{n-1}$ of the group $G$ with $|A_i|=q^i$, $i=1,\ldots,n-1$.
\end{example}
A matrix unit is a matrix with only one nonzero entry with value $1$.
The matrix unit with a $1$ in the $k$th row and $l$th column is denoted as $e_{kl}$. 
The rule of multiplication of matrix units: 
$$
e_{ij}e_{kl}=
\left\{
  \begin{array}{ll}
    e_{il}, & \hbox{if }j=k; \\
    0, & \hbox{otherwise},
  \end{array}
\right.
$$
where $0$ is the zero matrix.
The set $n\times n$ matrix units is a basis of the linear space of $n\times n$ matrices $M_{n,n}(\mathbb{F}_q)$.
If $U\subset M_{n,n}(\mathbb{F}_q)$, then by the symbol $\langle U\rangle$ 
we will denote the linear span of the set $U$.
Let 
$$
D_i=\langle e_{kl}\mid 1\leq k<l\leq n,\ l-k=i\rangle
$$
if $1\leq i\leq n-1$ and $D_i=\{0\}$ if $i\geq n$;
and
$$
V_i=D_{n-i}+\ldots+D_{n-1},\ i=1,\ldots,n-1.
$$
Let us denote the identity matrix with $e$.
We define for every $i\in\{1,2,\ldots,n-1\}$
$$
H_i=e+V_i=\{e+v\mid v\in V_i\}
$$ 
and 
$$
T_i=e+D_{n-i}=\{e+w\mid w\in D_{n-i}\}.
$$
It is not difficult to check that 
the sets $V_i$, $D_i$, $H_i$, and $T_i$ have the following properties ($1\leq i,j\leq n-1$):
\begin{enumerate}[label=$(\roman*)$]
  \item $\dim D_i=n-i$ and $\dim V_i=\binom{i+1}{2}$;
  \item $D_iD_j=D_{i+j}$;
  \item $V_iD_j=\left\{
                  \begin{array}{ll}
                    V_{i-j}, & \hbox{if $i>j$;} \\
                    \{0\}, & \hbox{otherwise;}
                  \end{array}
                \right.
        $
  \item  $V_iV_j=\left\{
                  \begin{array}{ll}
                    V_{i+j-n}, & \hbox{if $i+j>n$;} \\
                    \{0\}, & \hbox{otherwise.}
                  \end{array}
                \right.
        $ 
        \item $H_i$ is a normal subgroup of $G$;
        \item $H_1<H_2<\ldots<H_{n-1}=G$;
        \item $T_i$ is a right transversal to $H_{i-1}$ in $H_i$;
        \item $|T_i|=q^i$.
    \end{enumerate}

Now we construct the sets $A_i$.
Let $A_1=T_1$ and for each $i=2,\ldots,n-2$ let $h_i=e+e_{1,n-i}$ and
$$
A_i=T_ih_i=T_i+e_{1,n-i }.
$$
Further, it is important to remember that $n>3$. 
Let $t'=e+e_{1,n-1}$ and $s'=e+e_{12}+e_{23}$ and
$$
A_{n-1}=T_{n-1}\setminus\{e,e+e_{12}\}\cup\{t',s'\}.
$$
As an example, let us show what the sets of $A_i$ look like at $n=5$:
{\footnotesize
$$
A_1=
\left(
  \begin{array}{cccccc}
    1 & 0 & 0 & 0 & *  \\
      & 1 & 0 & 0 & 0  \\
      &   & 1 & 0 & 0  \\
      &   &   & 1 & 0  \\
      &   &   &   & 1  \\
  \end{array}
\right),
A_2=
\left(
  \begin{array}{cccccc}
    1 & 0 & 1 & * & 0  \\
      & 1 & 0 & 0 & *  \\
      &   & 1 & 0 & 0  \\
      &   &   & 1 & 0  \\
      &   &   &   & 1  \\
  \end{array}
\right),
A_3=
\left(
  \begin{array}{cccccc}
    1 & 1 & * & 0 & 0  \\
      & 1 & 0 & * & 0  \\
      &   & 1 & 0 & *  \\
      &   &   & 1 & 0  \\
      &   &   &   & 1  \\
  \end{array}
\right),
$$}
and
{\footnotesize
$$
A_4=\left\{
\left(
  \begin{array}{cccccc}
    1 & \alpha & 0 & 0 & 0  \\
      & 1 & * & 0 & 0  \\
      &   & 1 & * & 0  \\
      &   &   & 1 & *  \\
      &   &   &   & 1  \\
  \end{array}
\right),
\left(
  \begin{array}{cccccc}
    1 & 0 & 0 & 1 & 0  \\
      & 1 & 0 & 0 & 0  \\
      &   & 1 & 0 & 0  \\
      &   &   & 1 & 0  \\
      &   &   &   & 1  \\
  \end{array}
\right),
\left(
  \begin{array}{cccccc}
    1 & 1 & 1 & 0 & 0  \\
      & 1 & 0 & 0 & 0  \\
      &   & 1 & 0 & 0  \\
      &   &   & 1 & 0  \\
      &   &   &   & 1  \\
  \end{array}
\right)
\right\},
$$}
where for $\alpha$ in $A_4$ the following rule applies: if all asterisks are 0, then $\alpha\neq0,1$.

\section{Question}
As it is known, the converse of Lagrange's theorem is also true for supersolvable groups 
(see e.g. \cite[Lemma 4.1]{kabenyuk}).
Moreover, it is known that a finite supersolvable group is $k$-factorizable for every admissible integer $k$
\cite[Lemma 4.2]{kabenyuk}.
But in spite of this, the next question is open.

\begin{question}
    Is the complete factorization theorem valid for finite supersolvable groups?
\end{question}

Let us give some arguments in favor of answering this question in the affirmative.
We will say that subsets $A_1,\ldots,A_k$ of a group $G$ form a complete $(m_1,\ldots,m_k)$-factorization of $G$ if $A_1,\ldots,A_k$ form a complete factorization and $|A_i|=m_i$, $i=1,\ldots,k$.
An obvious consequence of proving the theorem is the following statement.
\begin{corollary}\label{corollary}  
    Let $G$ be a finite group and $|G|=m_1\ldots m_k$ where $m_1,\ldots,m_k$ are integers greater $1$ and $k\geq3$.
    If there exists a subnormal series 
    $$
    H_1<H_2<\ldots<H_k=G
    $$
    such that $|H_i|=m_1\ldots m_i$ for all $i=1,2,\ldots,k$,    
    then there exists a complete $(m_1,\ldots,m_k)$-factorization of group $G$.    
\end{corollary}

\begin{lemma}
If $p,q,r$ are distinct prime numbers and $G$ is a group of order $pqr$, 
then there exist a complete $(m_1,m_2,m_3)$-factorization of $G$, 
where $(m_1,m_2,m_3)$ is obtained from $(p,q,r)$ by an arbitrary permutation of coordinates.
\end{lemma}
\textit{Proof.}
Let, in fact, $P$, $Q$, $R$ be corresponding Sylow subgroups of $G$.
We can assume that $p>q>r$. 
In this case $P$ and $PQ$ are normal subgroups of $G$.
By Corollary \ref{corollary}, the group $G$ admits a complete $(p,q,r)$-factorization and then a complete $(r,q,p)$-factorization. 
Next, let $a_0\in P$, $b_0\in Q$, $c_0\in R$ be non-identity elements.
Since $P$ is a normal subgroup in $G$ we have
$$
P\cdot a_0b_0R\cdot Qc_0=a_0b_0PRQc_0=G \text{ and } G=a_0b_0R\cdot P\cdot Qc_0.
$$
Clearly, the sets $P$, $a_0b_0R$, and $Qc_0$ are pairwise disjoint, so that $G$ has complete 
$(p,r,q)$-, $(q,r,p)$-, $(q,p,r)$-, and $(r,p,q)$-factorizations.
The lemma is proved.


\begin{thebibliography}{9}
\bibitem[1]{Chin-Wang-Wong}
 A. Y. M. Chin, K. L. Wang and K. B. Wong,
\textit{Complete factorizations of finite abelian groups}, 
J. Algebra 628, 509 -- 523 (2023; Zbl 1516.20118)

\bibitem[2]{Hajos}
G.~\Hajos,
\textit{Uber einfache und mehrfache Bedeckung des n-dimensionalen Raumes mit einem Wurfelgitter},
Math. Z. 47, 427 -- 467 (1941; Zbl 0025.25401)

\bibitem[3]{kabenyuk}
 M.~Kabenyuk,
\textit{Factorizations of finite groups}, 2021,\\
\texttt{https://doi.org/10.48550/arXiv.2102.08605}

\bibitem[4]{Minkowski}
H.~Minkowski,
\textit{Geometrie der Zahlen},
Teubner, Leipzig,

\bibitem[5]{Rahimipour}
A.~R.~Rahimipour, A.~R.~Ashrafi, A.~Gholami,
\textit{The existence of minimal logarithmic signatures for some finite simple groups}.
Exp. Math. 27, 138 -- 146 (2018)

\bibitem[6]{Szabo}
S.~\Szabo,
\textit{Topics in Factorization of Abelian Groups}.
Birkhauser Verlag, Basel (2004)

\bibitem[7]{Szabo-Sands}
S.~\Szabo, A.~D.~Sands,
\textit{Factoring Groups into Subsets},
CRC Press, Taylor and Francis, New York, (2009; Zbl 1167.20030)
\end{thebibliography}
\end{document}